\documentclass[11pt]{amsart}
\usepackage{amssymb}
\usepackage{amsmath}
\usepackage{amscd}
\usepackage{amsfonts}
\usepackage{mathrsfs}
\usepackage{stmaryrd}
\usepackage{tikz}
\usetikzlibrary{arrows}
\usepackage[T1]{fontenc}
\usepackage{mathtools}
\usepackage{hyperref}
\usepackage{ulem}

\usepackage{bbm}

\newtheorem{prop-def}{Proposition-Definition}[section]

\begin{document}

\setlength{\oddsidemargin}{0cm} \setlength{\evensidemargin}{0cm}
\baselineskip=18pt

\begin{center}
{\bf{Systematic Counting of Restricted Partitions}}~\\
Mingjia YANG and Doron ZEILBERGER
\end{center}

\bigskip

{\bf{0.\quad Abstract}}~\\

We use  `partial difference operator schemes', and dynamical programming to design algorithms
that systematically count sets of integer partitions avoiding {\it any} set of patterns (of a certain,
natural, kind). We describe two approaches, a `negative' (adapting the Goulden-Jackson algorithm for enumerating words),
and a `positive' approach, that turns out to be much more efficient. Nevertheless the negative
approach has theoretical interest. ~\\

{\bf{1.\quad Introduction}}~\\

One of the cornerstones of enumerative combinatorics (and number theory!) are (integer) partitions. Recall that
a partition of a non-negative integer $n$ is a list of integers $(\lambda_1, \dots, \lambda_k)$ such that
$\lambda_1 \geq \dots \geq \lambda_k \geq 1$ and $\lambda_1+ \dots + \lambda_k=n$. 

As usuall, we will denote the number of integer partitions of $n$ by $p(n)$.
This is a very famous sequence,  OEIS sequence A41.

While there is no `explicit' formula for $p(n)$, there is a nice generating function, that goes back to Leonhard Euler.
Denoting the number of integer partitions of $n$ by $p(n)$, Euler discovered that
$$
\sum_{n=0}^{\infty} p(n)\,q^n \, = \, \prod_{i=1}^{\infty} \, \frac{1}{1-q^i} \quad .
$$

The {\bf bible} of the theory of partitions is George Andrews' classic [An]. We also strongly
recommend Drew Sills' fascinating monograph [S].

Suppose that you did not know about Euler's generating function, 
and you were given the task of computing the first, say, $1000$ terms of the sequence $p(n)$, how would you proceed?
The most straightforward way would be to try and use {\it dynamical programming}. Note that partitions have the {\it hereditary} property.
If you chop-off the largest entry of the partition of $n$, $(\lambda_1, \dots, \lambda_k)$, 
you would get a shorter partition, $(\lambda_2, \dots, \lambda_k)$, of $n-\lambda_1$.
Alas, because of the condition $\lambda_1 \geq \lambda_2$, we have to `remember' what $\lambda_1$ was, after kicking it out.
So we are {\bf forced} to consider a more general quantity, let's call it
$P(n,m)$, enumerating the set of partitions of $n$ whose largest part is {\bf exactly} $m$. Once we can compute this
more general quantity, the original object of interest, $p(n)$, is given by
$$
p(n) \, = \, \sum_{m=1}^{n} \, P(n,m) \quad .
$$

In order to compute $P(n,m)$ we have the obvious recurrence (alias partial {\bf difference equation})
$$
P(n,m) \, = \, \sum_{m'=1}^{m} P(n-m,m') \quad , \quad n \geq m \geq 1 ,
\eqno(FundamentalRecurrence)
$$
subject to the boundary conditions $P(m,m)=1$ and $P(n,m)=0$ if $n<m$. 
Replacing $n$ by $n-1$ and $m$ by $m-1$ in the above recurrence, and subtracting, one gets the even simpler recurrence
$$
P(n,m)= P(n-1,m-1)+ P(n-m,m) \quad .
\eqno(SimplifiedFundamentalRecurrence)
$$
This gives a {\bf quadratic time} (and quadratic memory) algorithm, $O(N^2)$ for compiling a table of $p(n)$ for $1 \leq n \leq N$.

This is not the most efficient way to compile such a table. An even better way is via Euler's recurrence (e.g. [An], p. 12)
$$
p(n)=\sum_{j=1}^{\infty} (-1)^{j-1} \left ( p(n-j(3j-1)/2)  + p(n-j(3j+1)/2) \right ) \quad,
$$
that was famously used by Major Percy MacMahon to compile such a table, that lead to Ramanujan's discovery of his
famous congruences (see [An]).

Already Euler considered the enumeration of sets of partitions obeying some {\it restrictions}. For example
the set of partitions into {\bf distinct} parts, let's call it $d(n)$, is given by the generating function ([An], p. 5)
$$
\sum_{n=0}^{\infty} d(n) q^n \, = \, \prod_{i=1}^{\infty} (1+q^i) \, = \,
 \prod_{i=0}^{\infty} \frac{1}{1-q^{2i+1}}
\quad .
$$

More recently, Rogers and Ramanujan (with the help of MacMahon, see [An] and [S]) considered the problem of enumerating
partitions with the property that the difference between consecutive parts is at least $2$, i.e. for which
$$
\lambda_i - \lambda_{i+1} \geq 2 \quad .
$$
The First Rogers-Ramanujan identity states that these numbers, let's call them $d_2(n)$, also have a nice
product generating function
$$
\sum_{n=0}^{\infty} d_2(n) q^n \, = \,  \prod_{i=0}^{\infty} \frac{1}{(1-q^{5i+1}) (1-q^{5i+4}) } \quad .
$$

We can say that distinct partitions {\bf avoid} the `pattern' $[a,a]$ and Rogers-Ramanujan partitions avoid both 
the pattern $[a,a]$ and the pattern $[a,a-1]$.

This naturally leads to the question of enumerating partitions avoiding an {\it arbitrary} (finite) set of patterns,
but first let's formally define the notion of a `pattern' in the context of partitions.

\bigskip

{\bf Definition.} A partition-pattern is a list $a=[a_1, \dots, a_r]$ of length $r \geq 1$ of {\it non-negative} integers.

{\bf Definition.} A partition $\lambda=(\lambda_1, \dots, \lambda_k)$ {\bf contains} the pattern
 $a=[a_1, \dots, a_r]$ if there exists $1 \leq i \leq k-r$ such that
$$
\lambda_{i}-\lambda_{i+1}=a_1 \quad,\quad
\lambda_{i+1}-\lambda_{i+2}=a_2 \quad,\quad
\dots
\lambda_{i+r-1}-\lambda_{i+r}=a_r \quad .
$$

For example, the partition $(7,6,5,4,4)$ contains the patterns $[1]$ (several times), the pattern $[0]$ (since $4-4=0$),
the pattern $[1,1]$ (because of $765$ and $654$), the pattern $[1,0]$ (because of $544$), the pattern $[1,1,1]$ 
(because of $7654$), the pattern $[1,1,0]$ (because of $6544$), and the pattern $[1,1,1,0]$.

{\bf Definition.} A partition $\lambda$ {\bf avoids} the pattern $a$ if it does {\bf not} contain the pattern $a$.

{\bf Definition.} A partition $\lambda$ {\bf avoids} the set of patterns $A$, if it avoids {\it every} pattern in $A$.

With this language, the class of distinct partitions are those that avoid the pattern $[0]$, while the
class of partitions whose differences are at least $2$ avoids the set of patterns $\{ [0],[1] \}$.

Our goal is to devise an efficient algorithm, that inputs an {\bf arbitrary} set of patterns, $P$,  and
an arbitrary positive integer $N$, and outputs the first $N$ terms of the sequence enumerating
partitions of $n$ avoiding the set of patterns $P$.

A natural approach is to adapt the celebrated Goulden-Jackson [GJ] method to this new context. Since it
is based on {\it sieving} (i.e. `signed-counting' using the deep identity $1+(-1)=0$) we call 
it a {\it negative} approach. 

The Goulden-Jackson method is lucidly explained (and significantly extended) in the article [NZ]. Recently it has been adapted [EZ] to counting
{\bf compositions} avoiding (a different kind of) patterns.

As it turned out, while this `negative' approach is very elegant, and of considerable {\it theoretical} interest,
it is less efficient than a more straightforward, `positive', approach, to be described later. In addition
to its theoretical value, the `negative' approach also serves as a good way to {\bf check} the correctness
of the far more efficient positive approach, in addition to using the `brute force'  of mere counting.
Since computer programs are still written by very unreliable human beings, it is always good to
have numerous checks, by making sure that the outputs to the {\it same} problem are always the same, using
several approaches, thereby empirically confirming all of them.

\bigskip

{\bf{2.\quad The ``negative'' approach}} \\

Recall that in the Goulden-Jackson Cluster method [GJ][NZ], 
one finds the weight enumerator for `marked words' 
and that turns out to be exactly the same as the target weight enumerator, that is, 
the weight enumerator for words avoiding a given set of subwords. 
Since the cluster method involves the {\it signed} counting of a larger set, and
often involves negative numbers,
we call it the "negative" approach here. 
However, in the setting of partitions, we cannot directly use the Goulden-Jackson Cluster method for the following reason: \\

In the Goulden-Jackson cluster method, one uses the important fact that if one peels off the first letter, or cluster,
of a (non-empty) marked word, then the result can be ANY marked word. So we have the following: \\

$M=\{empty\_word\} \, \cup \, V\,M  \, \cup \, C\, M$ 

(Note: $M$ is the set of all marked words, $V$ is the alphabet, $C$ is the set of all clusters) \\

Our basic idea is the same as in the Goulden-Jackson cluster method (we may call it the cluster method for simplicity from now on), however, since we are working with partitions, not words, we need the parts of the partition to be in non-increasing order. Therefore, when we peel off the first letter or cluster of a (non-empty) marked partition, the result is not any marked partition, but a marked partition with possibly a smaller first part such that after adding the cluster or the letter (that we peeled off) in front, it would still be a partition.\\

We also define weight a little differently than in the cluster method. Recall that in the cluster method, $weight(w,S)=(-1)^{|S|}s^{length(w)}$ ($S$ is the set of marks this word has). Here we define $weight(p,S)=(-1)^{|S|}s^{sum(p)}$ (where $sum(p)$ denotes the sum of the parts of $p$, that is, 
the integer that $p$ is partitioning.)\\

In order to use dynamical programming, we define the following: 

\begin{itemize}
    \item $P(A,k,m)$: the set of marked partitions that start with $k$ and having m parts, $A$ being the set of patterns to avoid.
    \item $C(A,k,l,w)$: the set of clusters starting with $k$, ending with $l$ and of width $w$, $A$ being the set of patterns to avoid. 
    \item $w(P(A,k,m))$: the weight enumerator of $P(A,k,m)$.
    \item $w(C(A,k,l,w))$: the weight enumerator of $C(A,k,l,w)$.

\end{itemize}
    
Let us start with a marked partition of largest part $k$ and $m$ parts. If the partition is empty ($m=0$), then the weight enumerator is $1$. If $m=1$, then the weight enumerator is $q^k$. If $m \geq 2$, the first part of the marked partition can be either part of a cluster or not, so for a fixed set of forbidden patterns $A$, we have the following decomposition: 

\begin{center} 
$P= kP\cup CP' $ 
\end{center} 

($P$ is the set of all marked partitions that start with $k$, having $m$ parts; $C$ is the set of all clusters starting with $k$, with width no greater than $m$; $P'$ is the set of marked partitions whose first part is no greater than the last part of clusters in $C$). More precisely, for $m \geq 2$, we have:

$$w(P(A,k,m)) = q^k \sum_{r=1}^{k}w(P(A,r,m-1))+\sum_{l=1}^{k}\sum_{w=1}^{m} (w(C(A,k,l,w)) \sum_{r=1}^{l} w(P(A,r,m-w))).$$


It remains to find $w(C(A,k,l,w))$. In order to do this, we introduce $w(C(A,v,k,l,w))$: the weight enumerator for clusters starting with $k$, with $v$ ($v \in A$) being the first pattern, ending with $l$ and of width $w$, $A$ being the set of patterns to avoid. For example, if $A=\{[2,1],[1,1]\}$, consider the cluster $\{8,6,5,3,2,1, \{[8,6,5],[5,3,2],[3,2,1]\} \}$. $v$ in this case would be $[2,1]$ (corresponding to the first mark $[8,6,5]$). It is apparent that $w(C(A,k,l,w))=\sum_{v \in A}w(C(A,v,k,l,w))$.\\

So how do we find $w(C(A,v,k,l,w))$? \\

For a given cluster, we have two scenarios: \\

($S1$) if the cluster has only one mark, then the weight for the cluster will just be $(-1) \cdot q^{sum(s)}$ ($s$ being the underlying partition). For example, the cluster $\{3,2,1, \{[3,2,1]\}\}$ has weight $-q^6$; \\

($S2$) if the cluster has more than one mark, we can peel off the first mark (leaving the overlapping part), and we get a smaller cluster. For example,  for the cluster $\{8,6,5,3,2,1, \{[8,6,5],[5,3,2],[3,2,1]\} \}$, after peeling off the first mark, we are left with the cluster $\{5,3,2,1, \{[5,3,2],[3,2,1]\}$. So, $weight(\{8,6,5,3,2,1, \{[8,6,5],[5,3,2],[3,2,1]\} \})=-q^{14}weight(\{5,3,2,1, \{[5,3,2],[3,2,1]\})$. \\

This is done in similar fashion as in the Goulden-Jackson cluster method. However, because of the nature of our extension, the details are more complicated. The first scenario occurs only if our input has width exactly 1 greater than the length of $v$, and the smallest part to ``match'' $k$ (the largest part) and the forbidden pattern, that is, $k=l+sum(v)$. To compute the weight for clusters in the second scenario, we first define $OVERLAP$, which takes two partitions $u$ and $v$ and outputs a set of lists. Each list is in the form $[q^i,j]$, where $j$ denotes the number of parts that $u$ and $v$ are overlapping, and $i$ denotes the sum of the parts of $u$ that is not overlapping with $v$. For example, $OVERLAP([4,3,2,2],[2,2,2,1])$ would return $\{[q^7,2],[q^9,1]\}$ because there are two possible ways of overlapping here. (Note: ``overlapping'' is defined in the usual sense, as in the cluster method, here the two possible overlaps are $[2,2]$ and $[2]$, the power 7 comes from $4+3$, the power 9 comes from $4+3+2$.) \\

Now, since we are really working with patterns (the $v$ in the input for $C(A,v,k,l,w)$ is a pattern, not a partition), we define $OVERLAP1$ which takes two patterns $u$ and $v$ and two integers $k1$ and $k2$ and let $u1$ and $u2$ be the corresponding partitions that start with $k1$ and $k2$ and with underlying pattern $u$ and $v$ respectively, and use $u1$ and $u2$ as input for $OVERLAP$.

For example, $OVERLAP1([1,1,0],[0,0,1],4,2)$ corresponds to $OVERLAP([4,3,2,2],[2,2,2,1])$ and also outputs $\{[q^7,2],[q^9,1]\}$. 

\newpage
Now we are ready to compute $w(C(A,v,k,l,w))$:\\

\begin{center}
$w(C(A,v,k,l,w))=(-1)q^{sum{\{v, k\}}}$(if  $k=l+sum(v)$ and $w=|v|+1$) $-\sum_{k1=1}^{k}\sum_{u \in A} \sum_{p \in OVERLAP1(v,u,k,k1)} p[1] \cdot w(C(A,u,k1,l,w-|v|-1+p[2]))$
\end{center}

\vspace{0.2in}
(Note: $\{v,k\}$ denotes the partition that start with $k$ and has underlying pattern $v$, for example, $\{[0,1],4\}=[4,4,3]$. $|v|$ is the length of the pattern $v$. $p[1]$ denotes the first part of the list $p$, $p[2]$ denotes the second part of $p$.)\\

In this formula, the part before the minus sign correspond to the first scenario, where the cluster have only one mark, and we will leave it to the reader to verify. If we are in the second scenario (computing the weight of the clusters that have more than one mark), we choose a pattern $u$ from $A$, and a largest part $k1$ ($1 \leq k1 \leq k$), and $\{u,k1\}$ is chosen to be the second mark of the cluster. We need to find all the ways $\{u,k1\}$ can overlap with $\{v,k\}$ (that is, compute $OVERLAP1(v,u,k,k1)$). Let us use the previous example $\{v,k\}=\{[1,1,0],4\}=[4,3,2,2]$, $\{u,k1\}=\{[0,0,1],2\}=[2,2,2,1]$. There are two ways they can overlap, if the overlap is $[2,2]$, then $p[1]$ would be $q^7$, and $p[2]$ would be 2. After chopping off the $[4,3]$ (that is, chopping off the first mark, leaving the overlapping part [2,2]) we would get a smaller cluster that starts with $k1=2$, still ends with $l$, and with width ($w-|v|-1+p[2])$, thus the formula above.\\

{\bf Remark:} One may wonder why we have to include the width as a variable. If we do not, and if [0] or [0,0], or [0,0,0] etc. is in $A$, then we would have infinitely many clusters (suppose there exist at least one cluster, we can then insert as many marks as we wanted in the middle) and we would have an infitely loop in our program. \\

{\bf{3.\quad The ``positive'' approach}} \\

While, for enumerating words (in a fixed alphabet) avoiding a given set of `patterns' (occurrences of consecutive subwords),
the negative approach, pioneered by Goulden and Jackson [GJ] is (usually) more efficient,
it turns out that this is not the case for the present problem of counting partitions avoiding the kind of patterns
discussed here.

The ``positive'' approach, to be described in this section, turns out to be much  more efficient than the negative approach
described in the previous section. Nevertheless, we believe that this partition analog of the Goulden-Jackson method is very elegant
and has theoretical interest. It is also possible that it may lead to more efficient approaches.

We use an extension of the dynamical programming approach described in the introduction that gave a quadratic-time and
quadratic memory algorithm to compute the original partition sequence $\{p(n)\}$, the iconic OEIS sequence $A41$.

It relied on the obvious fact that removing the largest part, $\lambda_1$, from a partition $\lambda=(\lambda_1, \dots, \lambda_k)$,
results in a smaller partition, $\lambda=(\lambda_2, \dots, \lambda_k)$, without extra conditions, except that $\lambda_2 \leq \lambda_1$.
That's why in the dynamical programming approach described in the introduction, we were {\bf forced} to compute
the more {\it refined} quantity, with {\bf two} arguments, $P(n,m)$, and that set-up the  {\it recurrence scheme} rolling.

If the set of forbidden patterns, $A$, consists only of patterns of length $1$,
$$
A= \{ [a_1], [a_2], \dots, [a_k] \} \quad ,
$$
then the analog of $(FundamentalRecurrence)$ is easy. Let $p_A(n)$ be the number of partitions of $n$ that avoid the patterns in the set $A$,
and let $P_A(n,m)$ be the number of such partitions whose largest part is $m$. Then

$$
P_A(n,m) \, = \, \sum_{
{
{1 \leq m' \leq m} \atop
{m-m' \not \in \{a_1, \dots, a_k\}}
}
}  P_A(n-m,m') \quad , \quad n \geq m \geq 1 \, ,
$$
and $p_A(n)=\sum_{m=1}^{n} P_A(n,m)$.

In order to motivate the general case, let's first do a simple special case, where we want to avoid
the single pattern $[1,1,1]$. In other words,
the set of patterns that we want to avoid is the singleton set $A=\{[1,1,1]\}$. 
Consider a  typical such partition $\lambda=(\lambda_1, \dots, \lambda_k)$,
whose largest part, $\lambda_1$, is $m$. If $\lambda_2 \neq \lambda_1 -1$, then removing $\lambda_1$ results
with the same type of partition, hence the number of partitions 
of $n$, that we are interested in, with $\lambda_1=m$ and $\lambda_2=m'$ is exactly the same as number
of such partitions of $n-m$ with largest part $\lambda_2$, since there is a one-to-one correspondence.
If you  have a good partition of $n-m$ with largest part $m'$, then sticking $m$ in the front can't cause trouble,
since $m-m' \neq 1$, so the forbidden pattern $[1,1,1]$ can't emerge.

On the other hand if $m'=m-1$ then we {\it can} create new trouble. If you have a partition of, $n$,  the form
$$
(m,m-1, \lambda_3, \dots, \lambda_k) \quad ,
$$
then the `be-headed' partition,of $n-m$
$$
(m-1, \lambda_3, \dots, \lambda_k) \quad ,
$$
must, {\bf in addition} to avoiding the pattern $[1,1,1]$ also avoid the pattern $[1,1]$ at the start.
This {\bf forces} us to introduce a new quantity, let's call it $P'_{[1,1,1]}(n,m)$ the number of partitions of $n$ with largest part
$m$,  avoiding the pattern $[1,1,1]$ {\it everywhere}, and in addition, avoiding the pattern $[1,1]$ at the very beginning
$$
P_{[1,1,1]}(n,m) \, = \, \sum_{
{{1 \leq m' \leq m}  \atop {m' \neq m-1}}
}  
P_{[1,1,1]}(n-m,m') \, + \, P'_{[1,1,1]}(n-m,m-1) \quad .
$$

We now need to set-up a scheme for $P'_{[1,1,1]}(n,m)$. If you have a partition of $n$ whose largest part
is $m$,  avoiding $[1,1,1]$, and in addition avoiding $[1,1]$ at the beginning, and the second largest part is $m'$ with
$m-m' \neq 1$, then removing the largest part, $m$, results in a partition of $n-m$ avoiding the pattern $[1,1,1]$,
and {\bf no conditions} at the beginning. On the other hand, if $m'=m-1$, then we have a partition of
$n-m$ with largest part $m-1$, avoiding $[1,1,1]$, and {\it in addition}, avoiding the pattern $[1]$ at the beginning. 
Let $P''_{[1,1,1]}(n,m)$ be the number of such partitions. We have

$$
P'_{[1,1,1]}(n,m) \, = \, \sum_{
{{1 \leq m' \leq m}  \atop {m' \neq m-1}}
}  
P_{[1,1,1]}(n-m,m') \, + \, P''_{[1,1,1]}(n-m,m-1) \quad .
$$

Similarly

$$
P''_{[1,1,1]}(n,m) \, = \, \sum_{
{{1 \leq m' \leq m}  \atop {m' \neq m-1}}
}  P_{[1,1,1]}(n-m,m') \, + \, P'''_{[1,1,1]}(n-m,m-1) \quad ,
$$
where  $P'''_{[1,1,1]}(n,m)$ is the number of partitions of $n$ with largest part $m$ avoiding the pattern $[1,1,1]$ and in addition
avoiding the empty list, $[]$, at the beginning. But this can never happen so $P'''{[1,1,1]}(n,m)$ is always zero. Note that we were forced to introduce two  auxiliary quantities, $P'(n,m)$, and $P''(n,m)$ that arose naturally.~\\

In general, for any given set of patterns $A$, the computer automatically sets-up a scheme, introducing
more general quantities, parameterized, in addition to the set of {\bf global} conditions $A$, by a set of
{\it local} conditions that should be avoided at the very beginning. Then, for each such set of
beginning restrictions, $A'$, depending on $m'$, either we are back to only the global conditions, $A$, i.e. the
new $A'$ is the empty set, or if $m-m'$ happens to be one of the starting entries of $A$ or $A'$, the
chopped partition, of $n-m$, in addition to obeying the global restrictions of $A$, must obey a brand-new
kind of restrictions $A''$. So each `state'  $(m,m',A')$ gives rise to a state $(m',m'',A'')$ for some
(possibly empty) set $A''$. Finding these ``children'' state is automatically done by the computer, setting
up a quadratic-time scheme. At the end of the day, we are only interested in the case where $A'=\emptyset$,
but we are forced to consider these auxiliary quantitities. Since there are only finitely many of them,
and there are still only two arguments (namely $n$ and $m$, where $1 \leq m \leq n$), 
the algorithm remains quadratic time and quadratic memory.
\newpage

\newpage
\bigskip
{\bf{4.\quad Maple packages}} \\

This article is accompanied by two Maple packages, {\tt RPneg.txt}, and {\tt RPpos.txt}, implementing the
two approaches described above. These are available from the front of this article

\url {http://sites.math.rutgers.edu/\~zeilberg/mamarim/mamarimhtml/rpr.html} \quad ,

where there is a large output file with the first $300$ terms of many sequences for many sets of forbidden patterns.
Most of them do not seem to be (as yet) in the OEIS, but some are.

The most important procedure in {\tt RPpos.txt} is $xnSeq(N,A)$, that outputs the
first $N$ terms in the enumerating sequence for partitions that avoid the set of patterns $A$.

$\bullet$ For $A=\{[0]\}$, we get, of course, the enumerating sequence for distinct partitions (alias odd partitions)
sequence $A9$ (\url{http://oeis.org/A000009}).

$\bullet$ For $A=\{[0],[1] \}$, we get the enumerating sequence for partitions whose minimal difference is at least $2$
(alias, via Rogers-Ramanujan, into parts $1,4$ modulo $5$), sequence $A3114$ (\url{http://oeis.org/A003114}).

$\bullet$ For $A=\{[1] \}$,  Sequence A116931, \url{http://oeis.org/A116931}, that goes back to MacMahon.

On the other hand the case $A=\{[2]\}$ is {\bf not} (yet) in the OEIS.

$\bullet$ For $A=\{[1], [0,0] \}$,  is sequence $A70047$, \url{http://oeis.org/A070047}.

On the other hand the case $A=\{[2], [0,0] \}$ is {\bf not} (yet) in the OEIS.

$\bullet$ The cases $A=\{[0,0] \}$,   $A=\{[0,0,0] \}$,   $A=\{[0,0,0,0] \}$,   are sequences $A726$, $A1935$ and $A35957$ respectively.

On the other hand the cases $A=\{[0,1]\}$, $A=\{[1,0]\}$ are not there (yet).

\bigskip

{\bf{5.\quad Future Work}} \\

The present algorithms assume that the parts can be anything. It would be 
fairly straightforward to impose congruence conditions on the parts.
A bit more challenging would be to have the restrictions also depend on congruences, for example,
Schur's celebrated 1926 theorem (see [An], p. 116), 
or the more complicated restrictions featuring in Shashank Kanade and Matthew C. Russell's intriguing 
conjectures ([KNR], [KR], [KR1], see also [S], pp. 149-152). This is an ongoing project.

\bigskip
{\bf{6. \quad Conclusion}} \\

We addressed the important problem of {\it systematically} enumerating classes of integer partitions 
avoiding {\it any} (finite) set of `patterns' (of a natural kind, that includes many classical cases).
We presented two algorithms, that we labeled `Negative' and `Positive'. The former is more elegant, while
the latter is more efficient.

\bigskip

{{\bf References}} \\

[An] George Andrews, {\it ``The Theory of Partitions''}, Addison-Wesley, Reading, MA, 1976. Reissued: Cambridge University Press, 1998. \\

[EZ] Shalosh B. Ekhad and Doron Zeilberger, {\it The ``Monkey Typing Shakespeare'' Problem for Compositions},
The Personal Journal of Shalosh B. Ekhad and Doron Zeilberger, Jan. 12, 2019. \hfill\break
{\tt http://sites.math.rutgers.edu/\~{}zeilberg/mamarim/mamarimhtml/kof.html} . \\

[GJ] Ian Goulden and David Jackson, {\it An inversion theorem for cluster decomposition of sequences with distinguished 
subsequences},   J. London Math. Soc.(2) {\bf 20} (1979), 567-576.  \\

[KNR]  Shashank Kanade, Debajyoti Nandi and Matthew C. Russell, {\it A variant of IdentityFinder and some new identities of Rogers-Ramanujan-MacMahon type}. {\tt https://arxiv.org/abs/1902.00790}. \\

[KR]  Shashank Kanade and  Matthew C. Russell, {\it IdentityFinder and some new identities of Rogers-Ramanujan type}, 
Experimental Mathematics {\bf 24} (2015),  419-423. 
\hfill\break
{\tt https://arxiv.org/abs/1411.5346} . \\

[KR1]  Shashank Kanade and  Matthew C. Russell, {\it Staircases to analytic sum-sides for many new integer partition identities of Rogers-Ramanujan
type}. Electron. J. Combin., 26(1):Paper 1.6, 2019. {\tt https://arxiv.org/pdf/1803.02515}  . \\

[NZ] John Noonan and Doron Zeilberger, {\it The Goulden-Jackson cluster method: extensions, applications, and implementations},
 J. Difference Eq. Appl. {\bf 5} (1999), 355-377.  \hfill\break
{\tt http://sites.math.rutgers.edu/\~{}zeilberg/mamarim/mamarimhtml/gj.html} . \\

[S] Andrew V. Sills, {\it ``An Invitation to the Rogers-Ramanujan Identities''}, CRC Press, Boca Raton, 2018. \\

\vfill\eject

\bigskip
\hrule

\bigskip

Mingjia Yang, Department of Mathematics, Rutgers University (New Brunswick), Hill Center-Busch Campus, 110 Frelinghuysen
Rd., Piscataway, NJ 08854-8019, USA. \hfill\break
Email: {\tt my237 at math dot rutgers dot edu}  . 

\bigskip

Doron Zeilberger, Department of Mathematics, Rutgers University (New Brunswick), Hill Center-Busch Campus, 110 Frelinghuysen
Rd., Piscataway, NJ 08854-8019, USA. \hfill\break
Email: {\tt DoronZeil at gmail  dot com}   .\\

\bigskip

{\bf Oct. 26, 2019}

\end{document}